\documentclass[12pt]{article}
\usepackage{amssymb,amsmath,cite}
\usepackage[dvips]{graphicx,color}
\usepackage{psfrag,subfigure}
\newtheorem{lemma}{Lemma}

\title{Analytic Continuation  of Divergent Integrals}
\author{Farhad Aghili\footnote{email: farhad.aghili@concordia.ca}}
\date{}

\begin{document}

\maketitle

\begin{abstract}
In this work, we investigate the improper integral of the monomial  \(
\mu(s) = \int_1^{\infty} x^{-s} \,dx \) as a continuous analogue of the infinite series representation of the Riemann $\zeta$-function, \(\zeta(s) = \sum_{n=1}^{\infty} n^{-s}\).
Both the monomial integral and the corresponding series converge for \(\mathrm{Re}(s) > 1\) and diverge for \(s \in \mathbb{C}\) with \(\mathrm{Re}(s) \leq 1\). In this paper, we construct an analytic continuation of the divergent monomial integral to the entire complex plane, excluding a simple pole at \(s = 1\), mirroring the analytic continuation of the $\zeta$-function. By performing term-by-term integration of the monomial over successive integer intervals and leveraging Newton’s generalization of the binomial theorem, we express the improper integral as a Dirichlet series. This approach establishes an elegant relationship between the \(\mu\)-function and the \(\zeta\)-function, leading to a functional equation that extends the divergent integral through analytic continuation and that the \(\mu\)-function is holomorphic everywhere except at \(s = 1\). 
\end{abstract}

\section{Introduction}

Analytic continuation is a powerful technique in mathematics used to assign meaningful finite values to infinite series that would otherwise diverge to infinity. In complex analysis, functions are initially defined on a small domain and then extended to a larger domain through analytic continuation. One of the most famous applications of this method is the Riemann zeta function. Riemann was the first to extend Euler's definition of the zeta function to the complex plane, establishing its analytic continuation to $\mathbb{C} \setminus \{1\}$:
\[
\zeta(s) = \sum_{n=1}^{\infty} n^{-s}.
\]
The analytic continuation of the zeta function has led to some of the most intriguing results in modern mathematics. For instance, the identity $\zeta(-1) = -\frac{1}{12}$ implies that the sum of all positive integers is somehow a finite number, which is neither positive nor an integer \cite{Berndt-1939}. Divergent series like this also appear in various calculations in modern physics, such as the Casimir force acting on two parallel plates due to vacuum energy. Remarkably, experimental measurements of the Casimir force align precisely with the result from the analytic continuation of the zeta function, specifically $\zeta(-3) = \frac{1}{120}$ \cite{Schumayer-2011, Elizalde-1995}. This astonishing agreement between theory and experiment demonstrates that analytic continuation is not just a mathematical abstraction but is deeply connected to the physical world. It seems that even nature “chooses” analytic continuation as a means to circumvent divergence.  

Similarly, divergent integrals also arise in modern physics, particularly in models describing interacting particles in quantum field theory \cite{Moreta-2013,Aghili-2019f,Schumayer-2011,Gobira-2021}. Techniques like regularization or renormalization are commonly used to extract finite values from divergent integrals by introducing cutoffs \cite{Felder-2016, Elizalde-Odintsov-1994}. Renormalization, in particular, is a fundamental tool in modern physics to remove the divergent part of integrals associated with physical quantities \cite{Salam-1951}. Other methods, such as weakly singular and hypersingular integral regularization, are based on the theory of distributions \cite{Zozulya-2007, Zozulya-2015}.

The aim of this paper is to extend the domain of certain divergent integrals to provide formal finite solutions through analytic continuation. Specifically, while the improper integral of the monomial function $x^{-s}$ converges only for $\mathrm{Re}(s) > 1$, the function
\[
\mu(s) = \int_1^{\infty} x^{-s} \, dx
\]
can be analytically continued and defined over $\mathbb{C} \setminus \{1\}$.

\section{Divergent Integral of Monomial}
The  improper integral $\int_1^{\infty} x^{-s} dx$ is convergent in domain $\{ s\in \mathbb{C}\; | \; \mbox{Re}(s) >1 \}$, i.e.,
\begin{equation}\label{eq:mu-def}
\int_1^{\infty} x^{-s} dx  =   \frac{1}{s-1}   \qquad   \mbox{Re}(s) >1
\end{equation}
Clearly, the improper integral of monomials becomes divergent if $\mbox{Re}(s) \leq 1$. Our end goal here is that the domain of the above improper integral as a function of complex variable $s$ is extended over the entire complex plane (except the isolated point $s=1$) in order to assign the divergent integral with a finite value in the sense of analytic continuation, i.e.,
\begin{equation} \label{eq:mu-sum}
\mu(s)=\int_1^{\infty} x^{-s} dx    \qquad s\in  \mathbb{C} \setminus \{ 1\}
\end{equation}
Notice that $\mu$-function definition \eqref{eq:mu-sum} can be viewed as the continuous analogy of the discrete Riemann zeta-function 
\begin{equation} \notag
\zeta(s) = \sum_{n=1}^{\infty} n^{-s}  \qquad s\in  \mathbb{C} \setminus \{ 1\}
\end{equation}
In the following analysis, we will fist show that there exists an analytical relationship between the $\mu$-function and $\zeta$-function.  Next, the link between the two functions is exploited to achieve analytic continuation of the $\mu$-function.  Let us the $\mu$-function represents the term by term integration of the monomial $x^{-s}$ over successive integer intervals, i.e.,
\begin{align} \notag
\mu(s) &= \int_1^{2} x^{-s} dx + \int_2^{3} x^{-s} dx + \int_3^{4} x^{-s} dx + \cdots \\ \label{eq:mu-Pn}
&= \sum_{n=2}^{\infty} P_n(s)
\end{align}
where $P_n(s)$ is the integral of the monomial over integer interval $n-1$ to $n$, i.e.,
\begin{equation} \label{eq:Delta1}
P_n(s) = \int_{n-1}^n x^{-s} dx = \frac{1}{1-s} \Big(  n^{1-s} - (n-1)^{1-s} \Big)
\end{equation}
According to the Newton's generalization of the binomial theorem, we can say
\begin{equation} \label{eq:n-1}
(n-1)^{1-s} =   \sum_{k=0}^{\infty}  (-1)^k {{1-s}\choose{k}} n^{1-s-k}    \qquad \forall s \in \mathbb{C}, \; n>1
\end{equation}
Here, the binomial coefficient is defined as
\begin{equation} \label{eq:binomial_coeff}
{{s}\choose{k}} = \frac{1}{k!}\prod_{j=0}^{k-1} (s-j)   = \frac{(s)_k}{k!}   \qquad s\in\mathbb{C}
\end{equation}
where $(\cdot)_k$ denotes falling factorial. Substituting \eqref{eq:n-1} in \eqref{eq:Delta1} yields
\begin{equation} \label{eq:Delta2}
P_n(s) = \frac{1}{s-1} \sum_{k=1}^{\infty}  (-1)^k  {{1-s}\choose{k}} n^{1-s-k}  \qquad \forall s \in \mathbb{C} \setminus \{ 1 \}, \; n\geq 2
\end{equation}
According to the ratio test, if the ratio of two successive terms of the series in \eqref{eq:Delta2} converges to $L$ which is less than one, i.e.,
\begin{equation} \label{eq:L}
L = \lim_{k \rightarrow \infty}\left| \frac{(1-s)_{k+1}}{n(k+1)(1-s)_k} \right| = \frac{1}{n}  \lim_{k \rightarrow \infty}\left| \frac{k+s-1}{k+1} \right| = \frac{1}{n},
\end{equation}
then the absolute convergence of the series is ensured. Clearly, \eqref{eq:L} implies that
\begin{equation} \notag
L < 1  \qquad \forall n\geq 2 \; s\in \mathbb{C}
\end{equation}
and thus the power series in \eqref{eq:Delta2}  converges absolutely with infinite radius of convergence.  Thus, $P_n(s)$ is always well-defined over the whole complex plane. 
Therefore, by virtue of \eqref{eq:mu-Pn} and \eqref{eq:Delta2}, the integral \eqref{eq:mu-def}  can be written as a Dirichlet series of the following form
\begin{align} \notag
\mu(s) = \int_1^{\infty} x^{-s} dx & = \frac{1}{s-1}   \sum_{k=1}^{\infty} (-1)^{k} {{1-s}\choose{k}}    \sum_{n=2}^{\infty} n^{1-s-k}  \\ \label{eq:int_sum}
&= \frac{1}{s-1}   \sum_{k=1}^{\infty} (-1)^{k} {{1-s}\choose{k}}  \big( \sum_{n=1}^{\infty} n^{1-s-k} -1 \big)
\end{align}
The infinite series in \eqref{eq:int_sum} diverges to infinity for $\mbox{Re}(s) \leq 1$. Nevertheless,  the Dirichlet series in the RHS of \eqref{eq:int_sum} is in the form of the Riemann zeta-function, i.e.,
\begin{equation} \label{eq:zeta}
\sum_{n=1}^{\infty} n^{1-s-k} = \zeta(k+s-1)
\end{equation}
Therefore, replacing  the equivalent of the Dirichlet from  \eqref{eq:zeta} into \eqref{eq:int_sum}, we arrive at the following expression of the improper integral in terms of the Riemann zeta-function
\begin{equation} \label{eq:mu}
\mu(s) =\frac{1}{s-1} \Big[   \sum_{k=1}^{\infty}(-1)^{k} {{1-s}\choose{k}}  \big(\zeta(k+s-1) -1 \big)      \Big],
\end{equation}
which is deemed to be new. Through the functional equation 
\begin{equation} \label{eq:functional}
\zeta(s) = \Big( 2^s \pi^{s-1} \sin \big( \frac{\pi s}{2} \big) \Gamma(1-s) \Big) \zeta(1-s)
\end{equation}
it has been proved that $\zeta(s)$  is holomorphic in the whole complex plane $\mathbb{C}$ except for
the only simple pole at $s = 1$. Thus, the meromorphic continuation of $\zeta(s)$ implies meromorphic continuation of $\mu(s)$. Now, we can proceed with simplification of the infinite series in  the RHS of \eqref{eq:mu}. For the sake of convenience definition, we split the $\mu$-function \eqref{eq:mu} in the following form
\begin{equation} \label{eq:mu_lambda}
\mu(s)  =\frac{1}{s-1} \lambda(s)
\end{equation}
where
\begin{equation}  \label{eq:lambda}
\lambda(s)=   \sum_{k=1}^{\infty}(-1)^{k} {{1-s}\choose{k}}  \big(\zeta(k+s-1) -1 \big)
\end{equation}
Since Riemann zeta function $\zeta(s)$ has a simple pole at $s=1$, the above expression implies that $\lambda(s)$ has poles at: $\{ 1 \} \cup \mathbb{Z}^- $, where  $ \mathbb{Z}^- =\{0,-1,-2,-3, \cdots \}$ denotes the set of non-positive integers. However, it will be later shown that those poles at  $\mathbb{Z}^-$ will be cancelled out by the zeros of the binomial coefficients. In the followings, we  will compute $\lambda(s)$ for  two complementary domains : i) $s\in \mathbb{C}-\big(\{ 1 \} \cup \mathbb{Z}^- \big)$, and ii) $s\in \mathbb{Z}^-$.

\begin{lemma}
The infinite sum of binomial coefficients and $\zeta$-function over $\mathbb{C} \setminus \big( \{ 1 \} \cup \mathbb{Z}^- \big)$ holds the following identity
\begin{equation} \label{eq:sum_zeta-1}
\lambda(s)=\sum_{k=1}^{\infty} (-1)^k {{1-s}\choose{k}}  \big( \zeta(k-s-1) -1 \big) =1  \qquad \forall  s\in \mathbb{C} \setminus \big( \{ 1 \} \cup \mathbb{Z}^- \big)
\end{equation}
\end{lemma}
{\sc Proof:} It was shown in  \cite{Shallit-1986} that Golback's theorem \cite{Titchmarsh-1951} assumes the following elegant form for the Riemann zeta function
\begin{equation} \label{eq:Golback}
\sum_{k=2}^{\infty} \big( \zeta(k) -1 \big) =1
\end{equation}
The alternative from of the aforementioned generalization of \eqref{eq:Golback} was presented in  \cite{Apostol-1954,Srivastava-1988}
\begin{equation} \label{eq:sum_zeta-2}
\sum_{k=1}^{\infty}{{k+s-2}\choose{k}} \big( \zeta(k+s-1) -1 \big) =1 \qquad \forall  s\in \mathbb{C} \setminus  \big( \{ 1 \} \cup \mathbb{Z}^- \big),
\end{equation}
which is not quite in the form of equation \eqref{eq:sum_zeta-1}. From the definition of binomial coefficient and falling factorial \eqref{eq:binomial_coeff}, we can say
\begin{align} \notag
(1-s)_k & = (1-s)(-s)(-s+1) \cdots (-s-k+2) \\ \notag
& = (-1)^k (s-1)s(s+1) \cdots (s+k-2) \\ \notag
& =  (-1)^k (s+k-2)_k
\end{align}
and hence
\begin{equation} \label{eq:2choose}
(-1)^k{{1-s}\choose{k}} =  {{k+s-2}\choose{k}}  \qquad \forall s \in \mathbb{C}
\end{equation}
By virtue of \eqref{eq:2choose} and \eqref{eq:sum_zeta-2}  we thus arrive  at the alternative form of the latter equation as presented in \eqref{eq:sum_zeta-1}.  $\Box$

In the reminder of this section, we compute $\lambda(s)$ for the second domain, i.e.,  $s\in \mathbb{Z}^-$.
\begin{lemma}
Sum of binomial coefficients and zeta-function hold the expression below:
\begin{equation} \label{eq:beta}
\beta(s)=\sum_{k=1}^{1-s}(-1)^{k} {{1-s}\choose{k}}  \big(\zeta(k+s-1) -1 \big) = \frac{(-1)^{-s}}{2-s} + 1 \qquad s\in\mathbb{Z}^-
\end{equation}
\end{lemma}
{\sc Proof:} One can verify that the successive binomial coefficients for positive integers $k$ and $-s$ hold the following useful identity
\begin{equation} \label{eq:successive_binomial}
{{1-s}\choose{k}} = \frac{2-k-s}{2-s} \cdot   {{2-s}\choose{2-k-s}}
\end{equation}
Moreover, the partial sums of the binomial coefficients satisfies the following identity \cite{Aupetit-2009}
\[ \sum_{j=0}^n (-1)^j {{n}\choose{j}} =0, \]
which implies
\begin{equation} \label{eq:sum_binomila=1}
\sum_{k=1}^{1-s}  (-1)^{-k} {{1-s}\choose{k}} = -1
\end{equation}
Upon substituting \eqref{eq:successive_binomial} and  \eqref{eq:sum_binomila=1} in the LHS of \eqref{eq:beta}, we can equivalently write the expression of the latter equation by
\begin{equation} \label{eq:beta_zeta}
\beta(s) = 1+ \frac{1}{2-s}\sum_{k=1}^{1-s} (-1)^{-k}{{2-s}\choose{2-k-s}}  (2-k-s) \cdot \zeta(k+s-1)
\end{equation}
Moreover, for positive integers $q \geq 0$, the Riemann zeta function is related to the Bernoulli numbers by
\begin{equation} \label{eq:zeta_benoulli}
\zeta(-q) = (-1)^q \frac{B_{q+1}}{q+1}
\end{equation}
Setting $q=1-k-s$ in \eqref{eq:zeta_benoulli} and then substituting the resultant equation in \eqref{eq:beta_zeta} yields the expression of the $\beta$ function as follows
\begin{equation} \notag
\beta(s) = 1 + \frac{(-1)^{1-s}}{2-s}  \sum_{k=1}^{1-s} {{2-s}\choose{2-k-s}} B_{2-k-s}
\end{equation}
or equivalently
\begin{equation} \label{eq:beta_bernoulli}
\beta(s) = 1 + \frac{(-1)^{1-s}}{2-s}  \sum_{k=1}^{1-s} {{2-s}\choose{k}} B_{k}
\end{equation}
On the other hand, the Bernoulli numbers satisfy the following property
\begin{equation} \notag
\sum_{k=0}^{q-1} {{q}\choose{k}} B_k =0            \qquad q \in \mathbb{Z},
\end{equation}
which implies
\begin{equation} \label{eq:sum_bernoulli}
\sum_{k=1}^{1-s} {{2-s}\choose{k}} B_k =-1
\end{equation}
Finally by plunging the equivalent value of the summation from \eqref{eq:sum_bernoulli} in the second term of the RHS of \eqref{eq:beta_bernoulli} one cane readily conclude identity \eqref{eq:beta}. $\Box$

\begin{lemma}
There exists the following relation between zeta function and the binomial coefficient
\begin{equation} \label{eq:alpha_lim}
\alpha(p) =\lim_{\epsilon \rightarrow 0} {{p-1+\epsilon}\choose{p}} \zeta(1-\epsilon) = -\frac{1}{p}  \qquad \forall p\in \mathbb{Z}^+
\end{equation}
\end{lemma}
{\sc Proof:}  Using the identities  $(r)_p=\Gamma(r+1)/\Gamma(r-p+1)$ and  $\Gamma(a+1)=a \Gamma(a)$ in the binomial coefficient  \eqref{eq:binomial_coeff}, one can expand the expression in \eqref{eq:alpha_lim} as follows
\begin{align} \notag
\alpha(p) &=  \lim_{\epsilon \rightarrow 0} \frac{ (p + \epsilon)_p \zeta(1- \epsilon) }{p !} \\ \notag
&= \lim_{\epsilon \rightarrow 0}  \frac{\Gamma(p+\epsilon) \zeta(1- \epsilon) }{p(p-1)! \Gamma(\epsilon)} \\ \notag
&= \lim_{\epsilon \rightarrow 0}  \frac{(p-1)! \zeta(1- \epsilon) }{p! \Gamma(\epsilon)} \\ \label{eq:zeta_gamma_ratio}
&= \frac{1}{p} \lim_{\epsilon \rightarrow 0}  \frac{ \zeta(1- \epsilon) }{\Gamma(\epsilon)}
\end{align}
Moreover, since the Reimann $\zeta$-function has a pole of first order at $s=1$, then it has a complex residue
\begin{equation} \label{eq:zeta_res}
\lim_{\epsilon \rightarrow 0} \epsilon \zeta(1+\epsilon) =1
\end{equation}
Similarly the Gamma function has a pole at $s=0$ and it is  well know that
\begin{equation} \label{eq:gamma_res}
\lim_{\epsilon \rightarrow 0} \epsilon \Gamma(\epsilon) =1
\end{equation}
Equations \eqref{eq:zeta_res} and \eqref{eq:gamma_res} indicate that zeta and Gamma functions approach equally fast towards their corresponding singular values at one and zero, respectively. Using \eqref{eq:zeta_res} and \eqref{eq:gamma_res} in \eqref{eq:zeta_gamma_ratio} concludes the proof of the Lemma, i.e., equation \eqref{eq:alpha_lim}.  $\Box$

Using the results of the above lemmas in \eqref{eq:beta_bernoulli} and \eqref{eq:alpha_lim}, and noting that 
\begin{equation}
{{1-s}\choose{k}} =0  \qquad  k> 2-s,
\end{equation}
we can now compute \(\lambda(s)\) for non-positive integers. This is achieved by decomposing the infinite sum in \eqref{eq:lambda} into three intervals: \( 1 \leq k \leq 1 - s \), \( k = 2 - s \), and \( k \geq 3 - s \), i.e.,
\begin{align} \notag
\lambda(s) &=  \sum_{k=1}^{\infty}(-1)^{k} {{1-s}\choose{k}}  \big(\zeta(k+s-1) -1 \big)  \qquad \forall s\in \mathbb{Z}^- \\ \notag
&= \underbrace{\sum_{k=1}^{1-s}(-1)^{k} {{1-s}\choose{k}}  \big(\zeta(k+s-1) -1 \big)}_{1 \leq k \leq 1-s } + \underbrace{(-1)^{-s}\alpha(2-s)}_{k=2-s} + \underbrace{0 + 0 + \cdots}_{k \geq 3-s} \\ \notag
& = \beta(s)+(-1)^{-s} \alpha(2-s)+ 0 \\ \label{eq:lambda_Z}
&= \Big( 1+ \frac{(-1)^{-s}}{2-s} \Big)  - \frac{(-1)^{-s}}{2-s} =1 \qquad \forall s\in \mathbb{Z}^-
\end{align}
We have now established all the necessary components to demonstrate that the \(\mu\)-function can be analytically extended to the entire complex plane, except for a simple pole at \( s = 1 \). From \eqref{eq:mu_lambda}, along with \eqref{eq:sum_zeta-1} and \eqref{eq:lambda_Z}, it follows that the analytic continuation of the divergent integral of the monomial \(\int_1^{\infty} x^{-s} \,dx\) is given by 
\begin{equation} \label{eq:def-mu}
\mu(s) \overset{\text{AC}}{=}  \frac{1}{s -1}  \qquad  \forall s\in \mathbb{C} \setminus \{ 1 \}.
\end{equation}

\section{Zero-Point Energy via Analytic Continuation}
In quantum field theory, the vacuum energy integral encapsulates the summation of all zero-point energy (ZPE) contributions within the vacuum state. However, this computation inherently yields divergent results, requiring the implementation of regularization and renormalization techniques to obtain physically meaningful, finite values by introducing cutoff scales \cite{Milonni-1993, Miransky-1994}. Although the conventional approach of introducing a cutoff parameter \( \Omega \) (e.g., a momentum or energy cutoff) can regulate the integral to yield a finite result \cite{Schwinger-1951, Elizalde-1994, MorenoPulido-2022}, the arbitrariness of \( \Omega \) often makes the outcome dependent on an artificial parameter rather than fundamental physical principles. As a case example, in this section, we will utilize the results of the analytic continuation of divergent integrals to calculate the zero-point energy (ZPE). The vacuum energy density for a  quantum massless  field  is typically given by the formula:
\begin{equation} \label{eq:ved_massless}
\rho = \frac{\hbar c}{4 \pi^2} \int_0^{\infty} k^3  \; d k,
\end{equation}
where $k$ is the wavevector of the quantum field, $\hbar =1.055 \times 10^{-34}$~kgm$^2$/s is the Planck reduced constant, and $c=3\times 10^8$~m/s is the speed of light. 
The integral in RHS of \eqref{eq:ved_massless}  diverges, indicating that the vacuum energy density for a massless field is formally infinite. In practical terms, a cutoff $\Omega$ is often introduced to replace the infinity as the upper bound of the integral in \eqref{eq:ved_massless} to limit the high-energy contributions 
\begin{equation} \label{eq:eps_massless}
\rho \approx \frac{\hbar c}{4 \pi^2} \int_0^{\Omega} k^3 dk.
\end{equation}
A high-energy cutoff is often set at the Planck scale, the energy level at which quantum gravitational effects are expected to become significant. The Planck length in term of Gravitational constant $G=6.674 \times 10^{-11}$~m$^3$/kg s$^2$ and other parameters are define   $l_p= \sqrt{\hbar G/c^3} = 1.6\times 10^{-35}$~m is related to the high-energy cutoff because it represents the smallest meaningful length scale in physics, and thus using the cutoff wavevector corresponding to the inverse of the Planck length, i.e., $\Omega =1/l_{p}$ could be used to calculate  the vacuum energy density \eqref{eq:ved_massless} with a finite value. Alternatively, we assume \eqref{eq:def-mu} to examine the derivation of the zero-point energy in \eqref{eq:ved_massless}. To achieve this, we apply the following variable substitution, shifting the integral's lower bound from zero to one to ensure consistency with the definition of the monomial integral \eqref{eq:mu-sum}
\begin{equation} \label{eq:change_var}
z^2 = 1+  k^2  \qquad  \mbox{with} \qquad z \; dz=k\; dk.
\end{equation}
That is
\begin{equation} \label{eq:eps_scaler}
\rho = \frac{\hbar c }{4 \pi^2}  \int_1^{\infty} (z^3 - z) dz   
\end{equation}
Finally, applying the analytic continuation of monomial integrals from \eqref{eq:def-mu} in \eqref{eq:eps_scaler} results in the following expression for the ZPE of massless fields, obtained without any cutoff:
\begin{equation} 
\rho = \frac{\hbar c }{4 \pi^2} \big( \mu(-3) - \mu(-1) \big) \overset{\text{AC}}{=} \frac{\hbar c }{16 \pi^2}. 
\end{equation}

\section{Conclusions}

In this paper, we applied the concept of analytic continuation to derive a unique finite solution for the divergent integral of the monomial across the entire complex plane, except at the pole $s = 1$. We treated the improper integral $\mu(s) = \int_1^{\infty} x^{-s} , dx$ as a continuous analogue of the infinite series representation of the Riemann zeta function, $\zeta(s) = \sum_{n=1}^{\infty} n^{-s}$, noting that both expressions diverge for $\mathrm{Re}(s) \leq 1$ but converge for $\mathrm{Re}(s) > 1$. We demonstrated that the monomial integral can be reformulated as a Dirichlet series through term-by-term integration over successive integer intervals, employing Newton’s generalization of the binomial theorem. This approach revealed an elegant relationship between the $\mu$-function and the $\zeta$-function, allowing us to assign a meaningful finite value to the divergent integral through analytic continuation.  As a case example, the analytic continuation of divergent monomial integrals has been used to evaluate zero-point energy as an alternative renormalization method.

\bibliographystyle{IEEEtran}

\bibliography{C:/Users/farha/OneDrive/Documents/Publications/bib/references}

\end{document}